\documentclass[12pt]{amsart}
\usepackage{amssymb}
\usepackage[mathcal]{eucal}
\usepackage[all]{xy}

\usepackage{amsfonts, amsmath, amscd}
\usepackage{enumerate}
\usepackage{pb-diagram}
\usepackage[usenames]{color}
\usepackage{graphicx}
\usepackage{amscd}

\newtheorem{deff}{Definition}[section]

\newtheorem{theorem}[deff]{Theorem}
\newtheorem{corollary}[deff]{Corollary}

\newtheorem{proposition}[deff]{Proposition}

\newtheorem{em-example}[deff]{Example}
\newtheorem{em-def}[deff]{Definition}        
\newtheorem{em-remark}[deff]{Remark}         
\newtheorem{em-question}[deff]{Question}

\newtheorem{problem}[deff]{Problem}
\newtheorem{definition}[deff]{Definition}

\newenvironment{example}{\begin{em-example} \em }{ \end{em-example}}

\newenvironment{remark}{\begin{em-remark} \em }{\end{em-remark}}

\newcommand{\EG}{EG_{k}}

\newcommand{\w}{\omega}

\def\span{\mathop{\rm span}}

\DeclareMathOperator*{\supp}{supp}

\DeclareMathSymbol{\res}{\mathord}{AMSa}{"16}

\def\:{\nobreak \hskip .1111em\mathpunct {}\nonscript \mkern
   -\thinmuskip {:}\hskip .3333emplus.0555em\relax}

\catcode`\@=12

\def\J{{\mathbb J}}

\def\N{{\mathbb N}}
\def\R{{\mathbb R}}
\def\Q{{\mathbb Q}}

\title[When is a locally convex space Eberlein-Grothendieck?]
{When is a locally convex space Eberlein-Grothendieck?}
\author{Jerzy K\c{a}kol, Arkady Leiderman}
\address{Faculty of Mathematics and Informatics, A. Mickiewicz University,
61-614 Pozna\'{n}, Poland and Institute of Mathematics Czech Academy of Sciences, Prague, Czech Republic}
\email{kakol@amu.edu.pl}

\address{Department of Mathematics, Ben-Gurion University of the Negev, Beer Sheva, P.O.B. 653, Israel}
\email{arkady@math.bgu.ac.il}

\keywords{Locally convex space, weak topology, $C_p(X)$ space, $C_k(X)$ space, compact space}
\subjclass[2010] {Primary 46A03; Secondary 46A20, 54C35, 54D30}
\date{\today}

\begin{document}
\begin{abstract}
In this paper we undertake a systematic study of those locally convex spaces $E$ such that $(E, w)$ is (linearly) Eberlein-Grothendieck,
where $w$ is the weak topology of $E$.

Let $C_{k}(X)$ be the space of continuous real-valued functions on a Tychonoff space $X$ endowed with the compact-open topology.
The main results of our paper are: (1) For a first-countable space $X$ (in particular, for a metrizable $X$)  the locally convex space
$(C_{k}(X), w)$ is Eberlein-Grothendieck if and only if $X$ is both $\sigma$-compact and locally compact;
 (2) $(C_{k}(X), w)$ is linearly Eberlein-Grothendieck if and only if $X$ is compact.

We characterize $E$ such that $(E, w)$ is linearly Eberlein-Grothen-\\
dieck
for several other important classes of locally convex spaces $E$. Also, we show that 
the class of $E$ for which $(E, w)$ is linearly Eberlein-Grothendieck preserves
 linear continuous quotients. Various illustrating examples are provided.
\end{abstract}
\thanks{The first named author
 is supported by the GA\v{C}R project 20-22230L and RVO: 67985840.}
\maketitle

\section{Introduction}\label{intro}
All topological spaces in the paper are assumed to be Tychonoff and all vector spaces are over the field
of real numbers $\R$.
By $C_k(X)$ and $C_p(X)$ we mean the space $C(X)$ of real-valued continuous functions defined on a Tychonoff space $X$
equipped with the compact-open and pointwise convergence topology, respectively.
For a locally convex space (lcs) $E$ we denote by $w$ the weak topology $w = \sigma(E, E^{\ast})$ of $E$ and, 
similarly, the $w^{\ast}$-topology of the dual $E^{\ast}$ is denoted by $w^{\ast}$.

The classic Eberlein-\v Smullian theorem states the equivalence of compactness, countable compactness and sequential compactness
for subsets $A\subset (E, w)$, where $E$ is any Banach space.
Remarkably, one can derive a proof of the Eberlein-\v Smullian theorem  from the analogous results about these 
versions of compactness in the function spaces $C_p(X)$.

Recall that a subset $A$ is {\it countably compact} in $Y$ if every countable sequence in $A$ has an accumulation point in $Y$.
The following significant statement which is intentionally formulated below in a simplified form is due to A. Grothendieck.

\begin{theorem}\label{Th_E-G} \cite[p. 107]{Arch_book} Let $X$ be a compact space, and let $A$ be a countably compact set in $C_p(X)$.
Then the closure of $A$ in $C_p(X)$ is compact.
\end{theorem}

Actually, Theorem \ref{Th_E-G} is valid assuming only that $X$ contains a dense $\sigma$-compact (i.e. a countable union of compact sets) subset. 
For further numerous developments of the Eberlein-Grothendieck theorem, we refer to \cite{Arch_book} and references therein.
Motivated by these facts A. V. Arkhangel'skii introduced the following definition. 

\begin{definition}\label{def:EG} \cite[p. 95]{Arch_book}
A topological space $Y$ is called an {\it Eberlein-Grothendieck space},
 if there is a homeomorphic embedding of $Y$ into the space $C_p(K)$ for some compact space $K$.
\end{definition}

Note here only that every metrizable space is Eberlein-Grothendieck \cite[Theorem IV.1.25]{Arch_book}.
Also, every Eberlein-Grothendieck topological space $Y$ has countable {\it tightness}, i.e. whenever a point $x$ is in the closure of a set $A \subset Y$,
then there is a countable $B \subset A$ such that $x$ is in the closure of $B$. Some questions about Eberlein-Grothendieck topological spaces
were considered in \cite{Aviles}.

In this paper we undertake a systematic study of those lcs $E$ such that $(E, w)$ is Eberlein-Grothendieck.
To the best of our knowledge this class of lcs has not been investigated before in such generality.
In particular, in Section \ref{section:1} we show that the class of lcs $E$ for which $(E, w)$ is Eberlein-Grothendieck
is closed under the operations of taking 1) linear subspaces; 2) linear continuous quotients; 3) countable products,  
but it does not preserve countable inductive limits.

Further, we observe that if $E$ is an Eberlein-Grothendieck lcs, then also $(E, w)$ is Eberlein-Grothendieck (Corollary \ref{cor_original}).
Example \ref{example_scattered} shows that the converse implication in general does not hold.

The following two fundamental results about the spaces $C_p(X)$ are due to O. Okunev \cite{Okunev}.
Recall that for $C_p(X)$ the weak and the original topologies are identical. 

\begin{theorem}\label{Th_Okunev} \cite{Okunev} If a lcs $H$ is a continuous open image of a subspace of $C_p(K)$ for some compact space $K$,
then the dual lcs $(H^{\ast}, w^{\ast})$ is $\sigma$-compact. 
\end{theorem}

\begin{corollary}\label{cor_Okunev} \cite{Okunev} A lcs $C_p(X)$ is Eberlein-Grothendieck if and only if $X$ is $\sigma$-compact.
\end{corollary}
We show that a similar assertion for the spaces $(C_k(X), w)$ fails.
Recall that a space $X$ is said to be {\it hemicompact} if there is a
sequence $\{K_n: n\in \omega\}$ of compact subsets of $X$ with the following property: if $K\subset X$ is compact then $K \subset K_n$ for some $n\in\omega$.
In Theorem \ref{Th1}, which is one of the main results of our work,
we prove that for every first-countable space $X$, the locally convex space  $(C_k(X), w)$ is Eberlein-Grothendieck if and only if $X$ is hemicompact.
As an easy corollary, we obtain a complete characterization of metrizable spaces $X$ such that
$(C_k(X), w)$ is Eberlein-Grothendieck (Theorem \ref{Th2}).

For a normed space $E$ there exists even a canonical {\it linear} embedding $T: (E, w) \to C_p(K)$,
 where $K$ is the closed unit ball in $E^{\ast}$ endowed with the weak$^{\ast}$-topology $w^{\ast}$.
 Being motivated by this simple observation we shall say that for a lcs $E$ the space $(E, w)$ 
 is a {\it linearly Eberlein-Grothendieck space} if $(E, w)$ can be linearly embedded into $C_p(K)$ for some compact space $K$.

In Section \ref{section:2} we show that for every $E$, a lcs $(E,w)$ is linearly Eberlein-Grothendieck if and only if the dual lcs $(E^{\ast}, w^{\ast})$ is compactly generated, i.e. 
$(E^{\ast}, w^{\ast})$ admits a compact subset $K$ whose linear span covers  $E^{\ast}$ (Theorem \ref{theor:dual_compgen}).
Establishing the relevant properties of compactly generated lcs we prove that 
the class of lcs $E$ for which $(E, w)$ is linearly Eberlein-Grothendieck, is closed
under the operation of taking linear continuous quotients (Theorem \ref{prop_operations_2}).

We characterize $E$ such that $(E, w)$ is linearly Eberlein-Grothendieck
for several important classes of lcs $E$. Along this line of research in Theorem \ref{th_C(X)} we prove that
a lcs $(C_k(X), w)$ is linearly Eberlein-Grothendieck if and only if $X$ is compact.

In the last Section \ref{section:3} we provide several illustrating examples. 
In particular, Example \ref{D(Omega)} demonstrates that both the space of test functions $\mathfrak{D}(\Omega)$
and the space of distributions $\mathfrak{D}'(\Omega)$ endowed with the weak topologies are not Eberlein-Grothendieck.

We recall definitions of several classes of locally convex spaces which will be used.
A {\it Fr\'echet space} is a completely metrizable locally convex space.
A lcs that is a locally convex inductive limit of a countable inductive system of Fr\'echet spaces is called a {\it (LF)-space}.
A lcs $E$ is said to be {\it Baire-like} 
if it is not the union of an increasing sequence of nowhere dense, balanced, convex sets \cite{Saxon}.
Clearly, every locally convex Baire space is Baire-like. 
For a lcs $E$ by the {\it strong dual} $E'$ of $E$ we mean the dual endowed with the strong topology $\beta(E', E)$. 

The dual space  $L_p(X) = C_p(X)^{\ast}$ plays a notable role in the paper. We use the fact that $X$ constitutes a closed Hamel basis of $L_p(X)$ (see \cite[Section 0.5]{Arch_book}).
Also, $L_p(X)$ is identical with the {\it free locally convex space} on $X$ endowed with the weak topology (see \cite{LMP}).

\section{Eberlein-Grothendieck lcs}\label{section:1}
The following characterization will be applied repeatedly in the sequel. 

\begin{proposition}\label{prop0} For every $E$, a lcs $(E, w)$ is Eberlein-Grothendieck if and only if the dual lcs $(E^{\ast}, w^{\ast})$ is $\sigma$-compact.
\end{proposition}
\begin{proof} $H = (E, w)$ is Eberlein-Grothendieck $\Longrightarrow$ $(E^{\ast}, w^{\ast})$ is $\sigma$-compact is an immediate consequence of Theorem \ref{Th_Okunev}.
The converse implication follows from the well-known fact that $(E, w)$ embeds into $C_p(L)$ by a linear homeomorphism, 
where $L$ is the space $(E^{\ast}, w^{\ast})$.
 Indeed, the required embedding is defined by the formula
 $$\xi: x \mapsto \xi_{x},\,\,\xi_{x}(x^{\ast})=x^{\ast}(x),\,\,x\in E, x^{\ast}\in E^{\ast}.$$
 Since $L$ is $\sigma$-compact we can apply Corollary \ref{cor_Okunev}.
\end{proof}

As a straightforward application of Theorem \ref{Th_Okunev} and Proposition \ref{prop0} we have

\begin{corollary}\label{cor_original} If $E$ is an Eberlein-Grothendieck lcs, then also $(E, w)$ is Eberlein-Grothendieck.
\end{corollary}

Later, in Example \ref{example_scattered} we show that the converse direction in Corollary \ref{cor_original} in general does not hold.
Since every metrizable topological space is Eberlein-Grothendieck we immediately obtain

\begin{corollary}\label{cor_metr} For every metrizable lcs $E$, the space $(E, w)$ is Eberlein-Grothendieck.
\end{corollary}

Another more informative argument implying Corollary \ref{cor_metr} will be presented later in the proof of Proposition \ref{EG2}.
The class of lcs $E$ with Eberlein-Grothendieck $(E, w)$  is invariant under certain basic topological operations.

\begin{proposition}\label{prop_operations_1} 
\mbox{}
\begin{itemize}
\item[{\rm (a)}] Let $(E, w)$ be an Eberlein-Grothendieck lcs. If there is a linear continuous quotient mapping $\pi$ from $E$ onto a lcs $F$, then
$(F, w)$ is also Eberlein-Grothendieck. 
\item[{\rm (b)}] Let $(E, w)$ be an Eberlein-Grothendieck lcs. Then $(F,w)$ is  Eberlein-Grothendieck
for every linear subspace $F \subset E$. 
\item[{\rm (c)}] Let $(E_n, w)$ be an Eberlein-Grothendieck lcs, where $n \in\omega$. Then the countable product $E = \prod_{n\in\omega} E_n$ also has the property that
$(E,w)$ is Eberlein-Grothendieck.
\end{itemize}
\end{proposition}
\begin{proof}
(a) We define the dual mapping $\pi^{\ast}: (F^{\ast}, w^{\ast}) \to (E^{\ast}, w^{\ast})$ by the formula: 
$$\pi^{\ast}(\phi) = \phi \circ \pi \,\, \mbox{for every}\,\, \phi \in F^{\ast}.$$
Since $\pi$ is a continuous and open surjection, the dual mapping $\pi^{\ast}$ establishes a linear homeomorphism between $(F^{\ast}, w^{\ast})$ and a closed subspace of $(E^{\ast}, w^{\ast})$
\cite[Theorems 8.12.1, 8.12.3]{Narici}.
The latter space is $\sigma$-compact, by Proposition \ref{prop0}, hence $(F^{\ast}, w^{\ast})$ is also $\sigma$-compact.
Finally, $(F, w)$ is Eberlein-Grothendieck, again by Proposition \ref{prop0}.\\
(b) It suffices to note that $(F, w)$ is a subspace of $(E, w)$ provided $F \subset E$ \cite[Theorem 8.12.2]{Narici}.\\
(c) If $E = \prod_{n\in\omega} E_n$ then $(E, w) = \prod_{n\in\omega} (E_n, w)$ (see \cite[Proposition 8.8.7]{Jarchow}).
Since $(E_n, w)$ homeomorphically embeds into $C_p(K_n)$ for some compact space $K_n$ for each $n \in \omega$, we have that
$(E, w)$ homeomorphically embeds into $C_p(X)$, where $X= \oplus_{n\in\omega} K_n$. It suffices to recall that the latter space $C_p(X)$ is Eberlein-Grothendieck
because $X$ is $\sigma$-compact.
The proof is complete.
\end{proof}

\begin{remark}\label{remark:phi}
Inductive limit of the countable sequence of lcs $(E_n)_n$, where each 
$(E_n, w)$ is Eberlein-Grothendieck, does not have to satisfy the same property.
Denote by $\varphi$ the $\aleph_0$-dimensional vector space endowed with the finest locally convex topology.
The space $\varphi$ can be identified with the strict inductive limit of the sequence of Euclidean spaces $\R^n$ \cite{Saxon}, and also 
$\varphi$ is isomorphic to the free locally convex space $L(\omega)$ on the countable discrete space $\omega$.
If $E = \varphi$ then the dual lcs $(E^{\ast}, w^{\ast})$ is isomorphic to the countable product $\R^{\omega}$ which clearly is not $\sigma$-compact.
Therefore, $\varphi$ endowed with the weak topology is not Eberlein-Grothendieck, by Proposition \ref{prop0}.
\end{remark}

For the brevity let us denote by $\EG$ the class of spaces $X$ such that a lcs $(C_k(X),w)$ is Eberlein-Grothendieck.
It is known that $C_k(X)$ is metrizable if and only if $X$ is hemicompact (see for instance \cite[Theorem 5.8.5]{Narici}).
Then, according to Corollary \ref{cor_metr}, $X \in \EG$ in the case that $X$ is hemicompact. 

\begin{problem}\label{problem1} Is the converse true, i.e. are the following properties equivalent:
(1) $(C_k(X),w)$ is an Eberlein-Grothendieck lcs and (2) $X$ is hemicompact?
\end{problem}

In this paper we show that the answer to Problem \ref{problem1} is positive for all  first-countable spaces $X$. 

For every $x\in X$ consider the evaluation mapping $\delta_x$,
 which is defined by the formula 
$$\delta_x(f)=f(x)\,\,\,\mbox{for every}\,\,\,f\in C(X).$$
It is well known that $\delta_x \in C_p(X)^{\ast} \subset C_k(X)^{\ast}$ for every $x \in X$.

We will extensively use the description of the elements from the dual space $C_k(X)^{\ast}$.
Below $\nu X$ denotes the Hewitt realcompactification of a Tychonoff space $X$. We identify the linear spaces $C(X)$ and $C(\nu X)$.
The essential of the following statement is well known and can be found in \cite{Schmets}.
Some relevant information can be found also in \cite{Velichko}.
The measures considered in the paper are {\it real-valued}.
The {\it support} $\supp(\mu)$ of a measure $\mu$ on $X$ is defined as the set of
all $x \in X$ such that every neighborhood $U$ of $x$ satisfies $|\mu|(U) > 0$.

\begin{theorem}\cite[Theorem III.3.3]{Schmets} \label{Th_Schmets}
Let $\phi$ be a continuous linear functional defined on $C_k(X)$.
Then there exists a unique Radon measure $\mu$ on $\nu X$ 
with the compact support $\supp(\mu) \subset X$ which fulfills the following properties
\begin{itemize}
\item[{(1)}] $\phi(f) = \int f d\mu$ for every $f\in C(X)$.
\item[{(2)}] If $f \in C(X)$ is such that $f = 0$ on $\supp(\mu)$, then $\phi(f)=0$. 
\item[{(3)}] For any open subset $U$ of $\nu X$ which satisfies the condition \\
$U \bigcap \supp(\mu) \neq \emptyset$, there exists some $f \in C(X)$ such that
						 $f = 0$ on $\nu X \setminus U$ and $\phi(f) = 1$.
\end{itemize}
\end{theorem}

\begin{proposition}\label{prop1} The subspace $\{\delta_x: x \in X\} \subset (C_k(X)^{\ast}, w^{\ast})$ is homeomorphic to $X$
and is closed in $(C_k(X)^{\ast}, w^{\ast})$.
\end{proposition}
\begin{proof} The first claim follows from the fact that the space $\{\delta_x: x \in X\}$ considered as a subspace of the double function space $C_p(C_p(X))$
is homeomorphic to $X$ (see \cite[Corollary 0.0.5]{Arch_book}).\\
The second claim is based on the description of the elements of $C_k(X)^{\ast}$ outlined in Theorem \ref{Th_Schmets} above. 
Let $\phi \in C_k(X)^{\ast}$ be defined by a measure $\mu$ on $\nu X$ with $\supp(\mu) \subset X$.
Assume that $\supp(\mu)$ contains at least two different points $x_1$ and $x_2$.
Then according to Theorem \ref{Th_Schmets} there are two continuous functions $f_1, f_2 \in C(X)$ such that 
$\supp(f_1) \bigcap \supp(f_2) = \emptyset$ and $\phi(f_1) = \phi(f_2) = 1$.
Define an open neighborhood of $\phi$ in $(C_k(X)^{\ast}, w^{\ast})$ as follows:
$$V= \{\zeta \in C_k(X)^{\ast}: \zeta(f_1) > 0, \zeta(f_2) > 0\}.$$ 
Then no functional $\delta_x$ belongs to $V$ since no $x\in X$ belongs both to $\supp(f_1)$ and $\supp(f_2)$.
Hence $\phi$ is not in the closure of $\{\delta_x: x \in X\}$. 
\end{proof}

As an immediate corollary we have
\begin{theorem}\label{cor1} If $X \in \EG$, then $X$ is $\sigma$-compact.
\end{theorem}
\begin{proof} By Proposition \ref{prop0} $(C_k(X)^{\ast}, w^{\ast})$ is $\sigma$-compact. It suffices to observe now that $X$ is homeomorphic to a closed subspace of $(C_k(X)^{\ast}, w^{\ast})$,
by Proposition \ref{prop1}.
\end{proof}

Proposition \ref{prop_operations_1} admits the following reformulation for the class $\EG$.
\begin{proposition}\label{prop_operations_2} 
\mbox{}
\begin{itemize}
\item[{\rm (a)}] Let $X \in \EG$. If $Y \subset X$ is closed, then also $Y \in \EG$.
\item[{\rm (b)}] Let $X \in \EG$. If $Y$ is an image of $X$ under a continuous compactly covering mapping, then also $Y \in \EG$.
\item[{\rm (c)}] Let $X_n \in \EG, n \in \omega$. Then the free topological sum $Y = \oplus_{n\in\omega}X_n$ also is in $\EG$.
\end{itemize}
\end{proposition}
\begin{proof}
(a) First, note that $X$ is Lindel\"of, by Theorem \ref{cor1}, hence $X$ is normal.
Our reasoning below is quite standard and is similar to analogous statements from \cite[Section 0.4]{Arch_book}.
Define the restriction mapping $\pi : C_k(X) \to C_k(Y)$, i.e. $\pi(f) = f\restriction_Y$ for all $f \in C_k(X)$.
Since $X$ is normal and $Y$ is closed in $X$, every function $g \in C_k(Y)$ can be extended to a function $f \in C_k(X)$.
Thus, $\pi$ is a mapping onto. The intersection of a compact subset $C \subset X$ with $Y$ is compact, therefore,
the mapping $\pi$ is continuous and open. Now Proposition \ref{prop_operations_1}(a) applies.\\
(b) Let $\gamma: X \to Y$ be a continuous compactly covering mapping.
Then the dual mapping $\gamma^{\ast}: C_k(Y) \to C_k(X)$ is a linear homeomorphic embedding, and Proposition \ref{prop_operations_1}(b) applies.\\
(c) Proposition \ref{prop_operations_1}(c) applies.
\end{proof}

Now we present one of the main results of our paper.

\begin{theorem}\label{Th1}
Let $X$ be a first-countable space. The following are equivalent
\begin{itemize}
\item[{\rm (a)}] $X$ is hemicompact.
\item[{\rm (b)}] $C_k(X)$ is metrizable.
\item[{\rm (c)}] $X \in \EG$, i.e. $(C_k(X), w)$ is Eberlein-Grothendieck.
\end{itemize}
\end{theorem}
\begin{proof} The equivalence (a) $\Longleftrightarrow$ (b) is very well known. (b) $\Longrightarrow$ (c) is true for every metrizable lcs.
By Corollary \ref{cor1} $X \in \EG$ $\Longrightarrow$ $X$ is $\sigma$-compact. Every $\sigma$-compact and locally compact space $X$ is hemicompact
(see for instance \cite{k_omega}); therefore, in order to prove (c) $\Longrightarrow$ (a) it suffices to show that $X \in \EG$ implies that $X$ is locally compact.
If $X \in \EG$ then $X$ is $\sigma$-compact, hence $X$ is paracompact.
A first-countable paracompact space is locally compact if and only if it does not contain
a closed subspace homeomorphic to the {\em metric fan} $M$ (see \cite[Lemma 8.3]{vD}).

The {\em metric fan} $M$ is a metrizable space defined as follows. As a set, $M$ is the union of countably many 
disjoint countable sequences $M_i  =\{x_{i n}:n \in \N\}, i \in \N$  plus a point $p$ "at infinity"; all points besides $p$ are isolated
in $M$, and a basic neighborhood $U_n$ of $p$ consists of $p$ and all points from $M_i$ such that $n \leq i$. 
Thus, the metric fan $M$ can be represented as a countable union of disjoint closed
discrete layers $M_i$, and a single non-isolated point $p$
such that for every choice $y_i \in M_i$ the sequence $(y_i)_i$ converges to $p$ in $M$.

On the contrary, assume that $X$ is not locally compact. Then $X$ contains a closed copy of $M$ and consequently $M \in \EG$ by Proposition \ref{prop_operations_2},
This is equivalent to the claim that $(C_k(M)^{\ast}, w^{\ast})$ is $\sigma$-compact by Proposition \ref{prop0}.
Our aim now is to show that the opposite is true: $(C_k(M)^{\ast}, w^{\ast})$ is not $\sigma$-compact because it
 contains a closed copy of the space of irrationals $\J \cong \N^\N$, which is not $\sigma$-compact.

By $\alpha \in \N^\N$ we understand a sequence of natural numbers $(\alpha(i))$. If $n \in \N$ then 
$\alpha\restriction_{n} = (\alpha(1), \alpha(2), \ldots, \alpha(n))$.
The space $M$ is countable, hence one can define a measure on $M$ by its value on each singleton.

For every element $\alpha \in \N^\N$ we define unequivocally the probability measure  $\mu_{\alpha}$ on $M$ by the following formulas:
$$\mu_{\alpha} (\{p\}) = 0\,\,\,\mbox{and}$$ 
$$
\mu_{\alpha} (\{x_{i n}\}) = \begin{cases}
\frac{1}{2^i},&\mbox{if}\,\,\,\alpha(i)=n,\\
0,&\mbox{otherwise}.
\end{cases}
$$

Each $\mu_{\alpha}$ defines a linear functional $\phi_{\alpha}$ on the linear space $C_k(M)$.
For every $\alpha \in \N^\N$  denote by $K_{\alpha}$ the converging sequence $(x_{i \alpha(i)}) \cup \{p\} \subset M$. 
Then  $K_{\alpha}$ is a compact support of $\mu_{\alpha}$ and $$|\phi_{\alpha}(f)| \leq \sup\{|f(x)| : x \in K_{\alpha}\}$$ for every $f \in C(X)$.
This means that $\phi_{\alpha}$ is a continuous linear functional on $C_k(M)$. Fix a notation $h$ for the correspondence $\alpha \mapsto \phi_{\alpha}$.
So we have a mapping $h$ from the space of irrationals $\J \cong \N^\N$ into the dual space $(C_k(M)^{\ast}, w^{\ast})$.
Denote the image of $h$ by $\Phi \subset (C_k(M)^{\ast}, w^{\ast})$.

We fix the notations for several test functions from $C(M)$ which will be used below.
For every $i \in \N$, a function $f_i \in C(M)$ is equal to 1 on the discrete layer $M_i$, and is zero for all other points.
If $x \in M \setminus \{p\}$, then $g_x \in C(M)$ is defined as follows: $g_x(x) = 1$ and $g_x(y) = 0$ for any other $y \in M$.\\
{\bf Claim 1.} $h$ is a 1-to-1 correspondence.\\
{\it Proof.} Let $\alpha \neq \beta$ and let $i$ be the first place such that $\alpha(i) \neq \beta(i)$.
Denote $x = x_{i \alpha(i)} \in  M \setminus \{p\}$.
Then $\phi_{\alpha} (g_x) \neq 0$, while $\phi_{\beta} (g_x) = 0$, meaning that $h(\alpha) \neq h(\beta)$.
\\
{\bf Claim 2.} $h$ is a homeomorphism between $\J $ and $\Phi$.\\
{\it Proof.} Pick any $\alpha \in \N^\N$. By definition of the product topology, the basic open neighborhoods of $\alpha$ in $\N^\N$ have the form
 $$U_n = \{\beta \in \N^\N: \beta\restriction_{n} = \alpha\restriction_{n}\},$$ where $n \in \N$.
For each $i \in \N$ denote by $y_i$ the point $x_{i \alpha(i)} \in M_i$.
 By our definition, $\phi_{\beta}(g_{y_i}) = \frac{1}{2^i}$ if $\beta(i) = \alpha(i)$ and $\phi_{\beta}(g_{y_i}) = 0$ if $\beta(i) \neq \alpha(i)$.
 For every $n\in\N$ we define an open neighborhood $V_n$ of $\phi_{\alpha}$ in
$\Phi$ by the following formula:
$$V_n = \{\phi_{\beta} \in \Phi: \phi_{\beta}(g_{y_i}) \neq 0, i =1, 2,\ldots, n\}.$$
It can be easily verified that $\phi_{\beta} \in V_n$ if and only if 
$\phi_{\beta}(g_{y_i}) = \phi_{\alpha}(g_{y_i}), i =1, 2,\ldots, n$,
if and only if $\beta\restriction_{n} = \alpha\restriction_{n}$. 

Now we show that the countable descending system of open neighborhoods $(V_n)_n$ produces a local base at $\phi_{\alpha}$ in $\Phi$.
It suffices to demonstrate that for every $f \in C(M)$ and every $\epsilon > 0$ there is $n\in \N$ such that
$V_n \subset W(f,\epsilon)=\{\phi_{\beta} \in \Phi: |\phi_{\alpha}(f) - \phi_{\beta}(f)| < \epsilon\}$.
Observe that a continuous function $f \in C(M)$ satisfies the following property: there are a real number $C > 0$ and a natural number $n$ such that
$|f(x)| < C$ for every $x \in \bigcup_{i > n} M_i$. Indeed, otherwise we would find an infinite set of indexes $I \subset \N$ and a sequence of points $(y_i)_i$ such that
$y_i$ belongs to the layer $M_i$ for every $i \in I$ with $|f(y_i)| \longrightarrow \infty$. This is not possible because the sequence $(y_i)_i$ converges to the point $p$ in $M$.
Without loss of generality we can assume also that $n$ is so big that $\frac{C}{2^{n-1}} < \epsilon$.
 Let $\phi_{\beta} \in V_n$. Then, as we have mentioned before, $\beta\restriction_{n} = \alpha\restriction_{n}$ and hence
$$|\phi_{\alpha}(f) - \phi_{\beta}(f)| = |\sum_{i=1}^{\infty}\frac{1}{2^i}(f(x_{i \alpha(i)}) - f(x_{i \beta(i)}))|=$$
 
$$=|\sum_{i=n+1}^{\infty}\frac{1}{2^i}(f(x_{i \alpha(i)}) - f(x_{i \beta(i)}))| \leq \sum_{i=n+1}^{\infty}\frac{1}{2^i} 2C = \frac{1}{2^{n-1}} C < \epsilon$$
This means that $\phi_{\beta} \in W(f,\epsilon)$. Thus, $V_n \subset W(f,\epsilon)$.
Finally we note that
$$h^{-1}(V_n) = U_n;\,\,h(U_n) = V_n\,\,\,\mbox{for every}\,\,\,n\in\N$$
and the proof of Claim 2 is complete. 
\\
{\bf Claim 3.} $\Phi$ is closed in $(C_k(M)^{\ast}, w^{\ast})$.\\
{\it Proof.} Pick any element $\xi \in C_k(M)^{\ast}$ and denote the support of a measure in $M$ which generates
the continuous linear functional $\xi$ by $K$.
First, if there is a layer $M_i$ such that the intersection $M_i \cap K$ is not precisely a singleton then $\xi$ is not in the closure of the set $\Phi$. 
Indeed, assume that $M_i \cap K = \emptyset$. Then for the function $f_i \in C(M)$ 
we obtain that $\xi(f_i) = 0$, while $\phi(f_i) = \frac{1}{2^i}$ for every $\phi \in \Phi$, meaning that $\xi$ can be separated from $\Phi$ by an open neighborhood inside $(C_k(M)^{\ast}, w^{\ast})$. 
If $M_i \cap K$ contains at least two different points, say $a$ and $b$, 
then we take two functions $g_a, g_b \in C(M)$. Put
 $$U = \{\zeta \in C_k(M)^{\ast}: \zeta(g_a) \neq 0, \zeta(g_b) \neq 0\}.$$
Then $U$ is an open neighborhood of $\xi$ in $(C_k(M)^{\ast}, w^{\ast})$, while $U \cap \Phi = \emptyset$.

Assume now that $\xi$ belongs to the closure of the set $\Phi$ in $(C_k(M)^{\ast}, w^{\ast})$. Denote the measure on $M$ which generates $\xi$ by $\nu$.
 It follows that the compact support $K$ of  $\nu$ meets each layer $M_i$ at precisely one point. 
 The latter fact means that $K$ is the converging sequence $(x_{i \alpha(i)}) \cup \{p\}$ for some unique element $\alpha \in \N^\N$.
 Our next step is to show that $\xi$ and $\phi_{\alpha} \in \Phi$ are identical. We look again at defined earlier functions $f_i \in C(M)$, for all $i \in \N$. 
Since $\phi(f_i) = \frac{1}{2^i}$ for every $\phi \in \Phi$, the same holds for $\xi$, i.e. $\xi(f_i) = \frac{1}{2^i}$.
By the structure of $\xi$ we infer that the measures $\mu_{\alpha}$ and $\nu$ are identical on each point from the layers $M_i$: 
$$
\nu(\{x_{i n}\}) = \begin{cases}
\frac{1}{2^i},&\mbox{if}\,\,\,\alpha(i)=n,\\
0,&\mbox{otherwise}.
\end{cases}
$$
It remains to ensure that also $\nu(\{p\}) = \mu_{\alpha}(\{p\}) = 0$. With this aim take the constant function $f = 1$. 
Observe that $\phi(f) =\sum_{i=1}^{\infty} \frac{1}{2^i} = 1$ for every $\phi \in \Phi$. Hence the same holds for $\xi$, i.e.  $\xi(f) = 1$.
But $\xi(f) = 1 + \nu(\{p\})$. Finally, $\nu(\{p\}) = 0$, measures $\mu_{\alpha}$ and $\nu$ are identical and the proof of Claim 3 is complete. 
\\
Summing up Claims 1 - 3 above, we conclude that $\Phi \cong \J$ is a closed subset of $(C_k(M)^{\ast}, w^{\ast})$. 
This provides a contradiction with our assumptions because $\J$ is not $\sigma$-compact. 
\end{proof}

As an easy corollary we obtain an exhaustive topological characterization of all metrizable spaces $X$ such that $(C_k(X), w)$ is Eberlein-Grothendieck.

\begin{theorem}\label{Th2}
Let $X$ be a metrizable space. The following are equivalent
\begin{itemize}
\item[{\rm (a)}] $X$ is hemicompact.
\item[{\rm (b)}] Either $X$ is compact, or there is a metrizable compact $K$ such that $X$ is homeomorphic to $K \setminus \{p\}$, where $p$ is a non-isolated point of $K$.
\item[{\rm (c)}] $C_k(X)$ is metrizable. 
\item[{\rm (d)}] $X \in \EG$, i.e. $(C_k(X), w)$ is Eberlein-Grothendieck.
\end{itemize}
\end{theorem}
\begin{proof}
In view of Theorem \ref{Th1} we need to show only (a) $\Longrightarrow$ (b). Every metrizable hemicompact space is locally compact (see \cite{k_omega}).
If $X$ is not compact, then a locally compact and metrizable space $X$ has a one-point metrizable compactification $K$.
\end{proof}

For countable metrizable $X$ we can represent $X \in \EG$ in a special form.

\begin{corollary}\label{cor_count} Let $X$ be a countable and first-countable space.
Then $X \in \EG$ if and only if there is a countable ordinal $\lambda$ such that either $X$ is 
homeomorphic to the closed segment of ordinals $[0, \lambda]$ or $X$ is 
homeomorphic to the interval $[0, \lambda)$.
\end{corollary}
\begin{proof} Every countable and first-countable space is metrizable; therefore, the statement follows from Theorem \ref{Th2} and
the known description of countable compact spaces.
\end{proof}

There is another important for applications class of lcs $E$ for which $(E, w)$ is Eberlein-Grothendieck if and only if $E$ is metrizable.
A lcs $E$ is called {\it Montel} if it is barrelled and every closed bounded set in $E$ is compact \cite[Section 11.5]{Jarchow}.
Recall that the strong dual of a Montel space is a Montel space.

\begin{proposition}\label{prop_Montel}
let $E$ be a Montel space. Then
$(E, w)$ is Eberlein-Grothendieck if and only if $E$ is metrizable.
\end{proposition} 
\begin{proof}
In light of Corollary \ref{cor_metr} we need to prove only one direction.
Assume that $(E, w)$ is an Eberlein-Grothendieck lcs. Denote by $E'$ the strong dual of $E$.
First, $E'$ is covered by countably many $w^{\ast}$-compact subsets, by Proposition \ref{prop0}. 
Second, since every Montel space is reflexive, the weak topology on $E'$ agrees with the $w^{\ast}$-topology.
Consequently, the barrelled space $E'$ is covered by a sequence of bounded sets.
 Using the argument from the proof of \cite[Proposition 2.12]{kak} we derive that  $E'$ has a fundamental sequence of bounded sets.
 Then the strong dual of $E'$ a metrizable lcs. 
However, the strong dual of $E'$ is equal to $E$ by reflexivity of $E$, and the proof is complete.
\end{proof}
\section{Linearly Eberlein-Grothendieck lcs}\label{section:2}
\begin{definition}\label{def:lin_EG}
A locally convex space $H$ is called {\it linearly Eberlein-Grothendieck} if 
 $H$ is linearly isomorphic to a subspace of the space $C_p(K)$ for some compact space $K$.
\end{definition}

It is easily seen that if a lcs $H$ satisfies Definition \ref{def:lin_EG} above, then the topology of $H$ coincides with its weak topology,
because this is true for every $C_p(K)$.
Recall also that $(E, w)$ is linearly Eberlein-Grothendieck for every normed space $E$. Hence any infinite-dimensional normed space $E$ provides an example
of a lcs which is not linearly Eberlein-Grothendieck while $(E, w)$ is linearly Eberlein-Grothendieck.

In this section we study lcs $E$ such that $(E, w)$ are linearly Eberlein-Grothendieck. 
Our first goal here is to find a counterpart of Proposition \ref{prop0} for this new stronger notion.

The following classic concept of topological algebra appears to be relevant to the subject.
A topological group $G$ which is algebraically generated by one of its compact subsets is called {\it compactly generated}.
By analogy, a topological vector space $L$ is called compactly generated if $L$ has a compact basis $K$, meaning that
the linear span of $K$ is equal to $L$. Evidently, a Banach space is compactly generated if and only if it is finite-dimensional.
 Nevertheless, the free locally convex space $L(X)$, as well as $L_p(X) = (L(X), w)$, over a compact space $X$ 
are compactly generated. We will show shortly that compactly generated lcs play an important role in our study.

\begin{theorem}\label{theor:dual_compgen} For every $E$, a lcs $(E,w)$ is linearly Eberlein-Grothendieck if and only if the dual lcs $(E^{\ast}, w^{\ast})$ is compactly generated.
\end{theorem}
\begin{proof} 
Assume first that there is an isomorphic embedding $\pi: (E, w) \to C_p(K)$, where $K$ is some compact space.
The dual of $C_p(K)$ is the space $L_p(K)$.
Then the adjoint mapping $\pi^{\ast}: L_p(K) \to (E^{\ast}, w^{\ast})$ is onto.
Therefore, the compact space $\pi^{\ast}(K)$ linearly generates $(E^{\ast}, w^{\ast})$.

Opposite direction. Assume that a compact subset $B$ linearly generates $(E^{\ast}, w^{\ast})$.
In the proof of Proposition \ref{prop0} we used a mapping of evaluation which always linearly embeds $(E, w)$ into $C_p(L)$, 
where $L$ is the whole space $(E^{\ast}, w^{\ast})$.
Since $B$ linearly generates $L$, here it is enough to take $C_p(B)$ as a target space.
Formally, the linear embedding $T: (E, w) \to C_p(B)$ is defined as follows:
$$(Tx)(x^{\ast})=x^{\ast}(x),\,\,x\in E, x^*\in B.$$
Below we explain why this embedding is a 1-to-1 mapping.
Take two different $x$ and $y$ in $E$. If we assume that $b(x) = b(y)$ for every $b \in B$,
then $x^{\ast}(x) = x^{\ast}(y)$ for every $x^{\ast} \in L$, which is of course false. So, the linear isomorphic embedding $T$
witnesses that $(E, w)$ is a linearly Eberlein-Grothendieck space.
\end{proof}

A closed subgroup of an arbitrary compactly generated abelian group in general does not have to be compactly generated.
Also, the discrete free group $F_2$ generated by two elements has a subgroup that is not finitely generated.
However, closed subgroups of a compactly generated locally compact abelian group are compactly generated (for the details see \cite{Ross}).

Does every closed linear subspace of a compactly generated lcs remain compactly generated?
Surprisingly, we were unable to find this question formulated and answered explicitly in any published source dealing with topological vector spaces.
By this reason we include its simple affirmative solution.

\begin{proposition}\label{prop_comp} Every closed linear subspace of a compactly generated lcs is also compactly generated.
\end{proposition}
\begin{proof} Assume that a lcs $L$ has a compact basis $K$ and  $H$ is a closed linear subspace of $L$.
Define $$K_n = \left\{\sum_{i=1}^n \lambda_i x_i: x_i \in K, \lambda_i \in [1, 1]\right\} \subset L.$$ 
Then each $K_n$ is a compact set because it is a continuous image of the product $[-1, 1]^n \times K^n$. 
 Denote by $H_n =  K_n \cap H$. Each $H_n$ is a compact subset of $H$ since $H$ is closed in $L$, and $\bigcup_{n\in\N} H_n = H$ since $\bigcup_{n\in\N} K_n = L$.
Finally, we define $$B = \bigcup_{n\in\N} \frac{1}{n^2} H_n \subset H.$$  
Clearly, $0 \in B$ and $\span(B) = H$ because $H_n  \subset \span(B)$  for every $n \in \N$.
It remains to show that $B$ is compact. Let $U$ be any absolutely convex neighborhood of zero in $L$. Every compact set is bounded, therefore, 
there is a natural number $n_0$ such that $K \subset n U$ for every $n > n_0$, or, equivalently, $\frac{1}{n} x \in U$ for every $x \in K, n > n_0$.
 Hence $\frac{1}{n^2} H_n \subset U$ for every $n > n_0$.
Indeed, every element of the set $\frac{1}{n^2} K_n$ is of the form $\sum_{i=1}^n \frac{\lambda_i}{n} y_i$, where each $\lambda_i \in [-1, 1], y_i \in U$
and $U$ is an absolutely convex set. We get that for every open cover of $B$ the neighborhood of zero covers 
all points of $B$ except perhaps for finitely many compact sets $\frac{1}{n^2} H_n$. Finite union of compact sets is covered by finitely many open sets 
and the proof is complete.  
\end{proof}

\begin{remark}\label{remark:free_tvs}
In general, Proposition \ref{prop_comp} is not true for non-locally convex spaces. Consider the free topological vector space $V(X)$, where 
$X$ is the closed unit interval $[0, 1]$. By construction, $V(X)$ is generated by the compact set $X=[0, 1]$. It is known that
the free topological vector space $V(Y)$ is isomorphic as a topological vector space to a closed vector subspace of $V([0, 1])$,
where $Y$ is the free union of countably many copies of the segment $[0, 1]$ (see \cite[Corollary 2.6]{LeiMor}).
The space $V(Y)$ is not compactly generated because every compact subset $K$ of $V(Y)$ is contained in a proper subspace $V(C) \subset V(Y)$
generated by a compact subset $C \subset Y$. Similarly, the free abelian group $A([0, 1])$ contains an isomorphic closed copy of $A(Y)$ which is not
compactly generated (see \cite{LMP}). This phenomenon cannot occur with the free locally convex space $L([0, 1])$, by the results of \cite{LMP}.
Note that the latter fact now is an immediate consequence of our Proposition \ref{prop_comp}.

It would be interesting to know whether the property established in Proposition \ref{prop_comp} characterizes compactly generated locally convex spaces among
 all compactly generated topological vector spaces.  
\end{remark}
 
\begin{theorem}\label{prop_operations_2} 
Let $(E, w)$ be a linearly Eberlein-Grothendieck lcs. If there is a linear continuous quotient mapping $\pi$ from $E$ onto a lcs $F$, then
$(F, w)$ is also linearly Eberlein-Grothendieck. 
\end{theorem}
\begin{proof}
As in the proof of Proposition \ref{prop_operations_1}, we define the dual mapping $\pi^{\ast}: (F^{\ast}, w^{\ast}) \to (E^{\ast}, w^{\ast})$ by the formula: 
$$\pi^{\ast}(\phi) = \phi \circ \pi \,\, \mbox{for every}\,\,  \phi \in F^{\ast}.$$
Since $\pi$ is a continuous and open surjection, the dual mapping $\pi^{\ast}$ establishes a linear homeomorphism between $(F^{\ast}, w^{\ast})$ and a closed linear subspace of $(E^{\ast}, w^{\ast})$.
The space $(E^{\ast}, w^{\ast})$ is compactly generated, by Theorem \ref{theor:dual_compgen}, hence $(F^{\ast}, w^{\ast})$ is also compactly generated, by Proposition \ref{prop_comp}.
Finally, $(F, w)$ is linearly Eberlein-Grothendieck, again in view of Proposition \ref{prop_comp}.\\
\end{proof}

The proof of the following statement is very similar to the proof of Proposition \ref{prop_operations_1} and is left to the reader.

\begin{proposition}\label{prop_operations_3}
\mbox{}
\begin{itemize}
\item[{\rm (a)}] Let $(E, w)$ be a linearly Eberlein-Grothendieck lcs. Then $(F,w)$ is linearly Eberlein-Grothendieck
for every linear subspace $F \subset E$. 
\item[{\rm (b)}] Let $(E_k, w)$ be an Eberlein-Grothendieck lcs, where $k= 1, 2, \dots, n$.
 Then the finite product $E = \prod_{k=1}^n E_k$ also has the property that
$(E, w)$ is linearly Eberlein-Grothendieck.
\end{itemize}
\end{proposition}

\begin{proposition}\label{prop_operations_4}
Let $E_k$ be lcs for all $k\in\omega$. 
Then for the countable product $E = \prod_{k\in\omega} E_k$ the space 
$(E, w)$ is not linearly Eberlein-Grothendieck.
\end{proposition}
\begin{proof} Countable product $E$ contains an isomorphic copy of $\R^{\omega}$.
The weak topology on $\R^{\omega}$ coincides with the original product topology. Hence, by Proposition \ref{prop_operations_3} (a), it suffices to observe that
$\R^{\omega}$ is not linearly Eberlein-Grothendieck. A proof of this fact can be extracted from the paper \cite{Uspen}.
A uniform space $F$ is called {\it precompact} if $F$ is a (uniform) subspace of a uniform compact space.
A uniform space which can be represented as a countable union of its precompact subspaces is called 
{\it $\sigma$-precompact}. V. Uspenskii \cite{Uspen} noticed that the property of being $\sigma$-precompact is preserved by any (uniform) subspace and $C_p(K)$ is
$\sigma$-precompact for every pseudocompact space $K$, while $\R^{\omega}$ is not $\sigma$-precompact. Hence, there is no a linear (even uniform) homeomoprhic embedding of 
$\R^{\omega}$ into $C_p(K)$ for a compact space $K$.

A somewhat alternative and direct argument is the following: $C_p(K)$ is covered by a sequence of closed bounded sets,
but $\R^{\omega}$ cannot be covered by a sequence of closed bounded sets, in view of the Baire Category Theorem.
\end{proof}

As an application of obtained  results, below we characterize $E$ such that $(E, w)$ are linearly Eberlein-Grothendieck
for several important classes of lcs $E$. 

\begin{proposition}\label{prop_Frechet}
Let $E$ be a Fr\'echet space. Then $(E, w)$ is linearly Eberlein-Grothendieck if and only if the metric of $E$ can be generated by a complete norm,
 i.e. $E$ is isomorphic to a Banach space.
\end{proposition}
\begin{proof}
 Every infinite-dimensional Fr\'echet space which is non-normable has the space $\R^{\omega}$ as a quotient space \cite{Eidelheit}.
Since $\R^{\omega}$ is not linearly Eberlein-Grothendieck, we deduce with the help of Theorem \ref{prop_operations_2} that a Fr\'echet space $E$ 
must be normable. Evidently, the generating norm is complete because the metric was supposed to be complete.
\end{proof}

In fact we can extend the last result to (LF)-spaces.

\begin{theorem}\label{th_LF}
Let $E$ be an (LF)-space. Then $(E, w)$ is linearly Eberlein-Grothendieck if and only if $E$ is normable.
\end{theorem}
\begin{proof}
We observe that $E$ does not contain a copy of $\varphi$ because $\varphi$ endowed with the weak topology is not Eberlein-Grothendieck,
by Remark \ref{remark:phi}. This allows us to apply \cite[Corollary 2.7]{kak} and to deduce that $E$ is a Baire-like lcs. 
By assumption, $(E, w)$ is isomorphic to a linear subspace of $C_p(K)$ for some compact space $K$. The space $C_p(K)$
is covered by countably many closed absolutely convex bounded sets $B_n = \{f(x) \in C(X): \sup_{x\in K} |f(x)| \leq n\}$, where $n\in \omega$.
Hence, $(E, w)$ is also covered by such a sequence of bounded sets $B_n \cap E$.
The sets $B_n$ remain closed absolutely convex and bounded in the original topology of $E$.
Since $E$ is Baire-like, we conclude that $E$ admits a bounded neighborhood of zero. This implies that $E$ is normable.
\end{proof}

\begin{remark}\label{rem_LF} An (LF)-space is Baire-like if and only if it is metrizable \cite[Corollary 2.7]{kak}. However, 
a normed (LF)-space does not have to be complete (see \cite[Section 8.7]{Bonet}).
\end{remark}

\begin{theorem}\label{th_C(X)}
For any Tychonoff space $X$ the following conditions are equivalent:
\begin{itemize}
\item [(i)] $X$ is compact;
\item [(ii)] $C_p(X)$ is linearly Eberlein-Grothendieck;
\item [(iii)] $(C_k(X), w)$ is linearly Eberlein-Grothendieck.
\end{itemize}
\end{theorem}
\begin{proof} We show only non-trivial implications (ii) $\Longrightarrow$ (i); (iii) $\Longrightarrow$ (i).
Assuming (ii) we know that $X$ must be $\sigma$-compact, by Corollary \ref{cor_Okunev}. On the other hand, $X$ must be pseudocompact because otherwise $C_p(X)$ would contain
an isomorphic copy of $\R^{\omega}$ \cite{Uspen}, which is not linearly Eberlein-Grothendieck. Every $\sigma$-compact and pseudocompact space is compact and the proof of (i) is finished.

Similarly, assuming (iii) we know that $X$ must be $\sigma$-compact, by Theorem \ref{cor1}. Additionally,
$X$ must be pseudocompact because otherwise $C_k(X)$ would contain
an isomorphic copy of $\R^{\omega}$, this time by \cite[Theorem 2.12]{kak} and then $(C_k(X), w)$ would contain
an isomorphic copy of $\R^{\omega}$ as well.
As before, this means that $X$ is compact and the proof of (i) is finished.
\end{proof}

Note that while $(E, w)$ is Eberlein-Grothendieck for every metrizable $E$, we have plenty of examples of 
metrizable $E$ such that $(E, w)$ is not linearly Eberlein-Grothendieck.
 Nevertheless, we have  the following

\begin{proposition} \label{EG2}
If $E$ is a metrizable lcs, then
 $(E, w)$ embeds linearly into $C_p(K, \R^{\omega}) \cong C_p(X)$, where $K$ is a compact space and
$X$ is the free sum of countably many copies of $K$.
\end{proposition}
\begin{proof}
Let $(U_{n})_{n}$ be a decreasing  base of absolutely convex and closed neighborhoods of zero in $E$.
For each $n\in\omega$ define the polar set 
$$U_{n}^{\circ}=\{x^{\ast}\in E^{\ast}: |x^{\ast}(x)|\leq 1, x\in U_{n}\}.$$ 
Then  each $U_{n}^{\circ}$ is a $w^{\ast}$-compact set in $E^{\ast}$ by the classic Alaoglu-Bourbaki theorem, and $E^{\ast}=\bigcup_{n\in\omega}U_{n}^{\circ}$.
Hence  $K=\prod_{n\in\omega}U_{n}^{\circ}$ is a compact set in the space $\prod_{n\in\omega}(E^{\ast}_{n})_{w{\ast}}$, where $E^{\ast}_{n}=E^{\ast}$ for each $n\in\omega$. 
Define a linear mapping $\xi: x \mapsto \xi_{x}\in C_p(K, \R^{\omega})$, $x\in E$, by the formula
$$\xi_{x}(x^{\ast})=(x^{\ast}_{n}(x)),\, \mbox{where}\,\, x^{\ast}=(x^{\ast}_{n})\in \prod_{n\in\omega}U^{\circ}_{n}.$$

We claim first that $\xi$ is \emph{injective}. Indeed, assume that $\xi_{x}=\bar{0}$ for $x\in E$, where $\bar{0}$ denotes the zero of $C(K, \R^{\omega})$.
 Take any $x^{\ast}\in E^{\ast}$.
 Then there exists $n\in\omega$ such that $x^{\ast}\in U^{\circ}_{n}$.
  Set $y^{\ast}=(0, ..., 0, x^{\ast}, 0, ...)\in U^{\circ}_{n}\times\prod_{m\neq n}U^{\circ}_{m}$.
Since $\xi_{x}=\bar{0}$, we  know  that  $\xi_{x}(y^{\ast})=0$.
Hence $(0, ..., x^{\ast}(x), 0, ...) = 0$, which implies that  $x^{\ast}(x)=0$ for each $x^{\ast}\in E^{\ast}$, meaning that  $x$ is zero in $E$, so $\xi$ is injective.

Now we prove that $\xi$ is \emph{continuous.} If $(x_{\gamma})$ is a net in $E$ which converges to zero in $(E, w)$, then $x^{\ast}(x_{\gamma})\rightarrow 0$ for each $x^{\ast}\in E^{\ast}$. 
 Let $y^{\ast}=(x^{\ast}_{n})\in\prod_{n\in\omega}U^{\circ}_{n}$. Since $x^{\ast}_{n}(x_{\gamma})\rightarrow 0$ for each $n\in\omega$,
 then $y^{\ast}(x_{\gamma})\rightarrow 0$, showing that $\xi_{x_{\gamma}}\rightarrow \bar{0}$.

In order to show that the inverse to $\xi$ is continuous, assume that $\xi_{x_{\gamma}}\rightarrow \bar{0}$. We need to show that $x_{\gamma}\rightarrow 0$ in $(E, w)$. Take any $x^\ast\in E^\ast$.
 As above, there exists $n\in\omega$  such that $y^{\ast}=(0, ..., 0, x^{\ast}, 0, ...)\in U^{\circ}_{n}\times\prod_{m\neq n}U^{\circ}_{m}$.
   By assumption we have that $\xi_{x_{\gamma}}(y^{\ast})\rightarrow 0$.
	This implies that  $\xi_{x_{\gamma}}(y^{\ast})=(0, ..., x^{\ast}(x_{\gamma}), 0, ...)\rightarrow 0,\,\text{hence}\,\, x^{\ast}(x_{\gamma})\rightarrow 0$.

The last claim that $C_p(K, \R^{\omega})$ is canonically isomorphic with $C_p(X)$ is well known (see \cite[Propositions 0.3.3 and 0.3.4]{Arch_book}).
\end{proof} 

\section{Illustrating examples}\label{section:3}
\begin{example}\label{example_scattered}
Let $X$ be any uncountable scattered compact space. By $C(X)$ here we mean a Banach space with the $\sup$-norm.
We define a lcs $E$ to be $C(X)$ endowed with the coarser vector topology $\tau$ of uniform convergence on the countable compacts subsets of $X$.
We claim that the spaces $(E, w)$ and $(C(X), w)$ coincide.
Indeed, every linear continuous functional $\xi$ on $C(X)$ for a scattered compact space $X$ is defined by an atomic measure (which has a countable and compact support) on $X$
(see \cite[Corollary 19.7.7]{Semadeni}). This means that $\xi$ remains continuous being considered as a linear functional on $E$.  
So, $E^\ast = C(X)^{\ast} = \ell_1(X)$. Since $(E, w)$ and $(C(X), w)$ coincide and $(C(X), w)$ obviously is linearly Eberlein-Grothendieck,
we conclude that $(E, w)$ is also linearly Eberlein-Grothendieck. 
However, $E$ with the original topology $\tau$ does not have to be Eberlein-Grothendieck. 
Take any uncountable scattered compact space $X$ with the following properties: 1) There is an uncountable subset $Y \subset X$ consisting of isolated points;
2) Every separable closed subset of $X$ is at most countable.
 For instance, let $X$ be the compact ordered space of ordinals $[0, \omega_1]$,
or let $X$ be the one-point compactification of an uncountable discrete set $Y$. For each $y\in Y$ define a function $f_y \in C(X)$ as follows:
$f_y(x) = 1$ if $x=y$ and $f_y(x) = 0$ otherwise. One can easily verify that the constant zero function belongs to the closure of 
$\{f_y: y\in Y\}$ in $E$ but zero does not belong to the closure of 
$\{f_y: y\in A\}$ in $E$ for every countable $A\subset Y$.
It appears that $E$ does not have countable tightness, hence $E$ cannot be homeomorphically embedded into $C_p(K)$ for any compact space $K$.
\end{example}

\begin{example}\label{example_count}
The implication $X \in \EG \Longrightarrow X$ is locally compact in Theorem \ref{Th1} and
Corollary \ref{cor_count} are not valid without the assumption that $X$ is first-countable.
Let $p$ be any point from the remainder of the Stone-\v{C}ech compactification $\beta{\N}\setminus \N$.
Put $X= \N \cup \{p\} \subset \beta{\N}$. Then $X$ does not contain any infinite compact subset. Therefore, $C_k(X)$ coincides with the metrizable lcs $C_p(X)$.
 However, $X$ is not locally compact.
So, the general Problem \ref{problem1} appeals for a complete answer even for countable spaces $X$.
\end{example}

\begin{example}\label{example_restr)} Let $X$ be the compact ordered space of ordinals $[0, \omega_1]$ and $Y$ be its subspace $[0, \omega_1)$.
Then $X$ is compact, while $Y$ is not compact and every continuous real-valued function defined on $Y$ extends continuously to $X$.
These very well-known facts show that the linear mapping of restriction $\pi$ maps $C_p(X)$ onto $C_p(Y)$, and 
$(C_p(Y), w)$ is not linearly Eberlein-Grothendieck despite that $(C_p(X), w)$ is linearly Eberlein-Grothendieck, by Theorem \ref{th_C(X)}.
We conclude that Theorem \ref{prop_operations_2} is not valid in general if the mapping $\pi$ is not open.
\end{example}

\begin{example}\label{example_sec-count}
Let $X$ be any second-countable non-locally compact space, for example, the space of rationals $\Q$.
Then $(C_k(X), w)$ is not Eberlein-Grothendieck by Theorem \ref{Th2}. Observe that
$C_k(X)$ is a so-called $\aleph_0$-space, therefore, $(C_k(X), w)$ has countable network as a continuous image of $C_k(X)$, and hence 
$(C_k(X), w)$ has countable tightness
(for the details see \cite{Michael}).
\end{example}

\begin{example}\label{example_L(X)}
A lcs $X$ is defined to have a {\em neighborhood $\w^\w$-base} at zero if there exists
a neighborhood base $\{U_\alpha: \alpha\in\w^\w\}$ at zero such that $U_\beta\subset U_\alpha$ for all
elements $\alpha\le\beta$ in $\w^\w$.
 Every metrizable lcs has an $\omega^\omega$-base at zero and the properties of lcs possessing an $\omega^\omega$-base at zero
resemble the properties of metrizable LCS \cite{kak}. However, Corollary \ref{cor_metr} cannot be generalized to $\omega^\omega$-based lcs.
Let $E$ be be the free locally convex space $L(X)$. It has been known that $L(X)$ has a $\omega^\omega$-base at zero
for every metrizable compact space $X$ (see \cite{GKL}). But $L(X)$ equipped with its weak topology, $L_p(X)$, is an Eberlein-Grothendieck space
if and only if $X$ is finite, by \cite[Corollary 1]{Okunev}. In particular, for the countable discrete space $X = \omega$ we have that
$L(X) \cong \varphi$, and $(L(X)^{\ast}, w^{\ast}) \cong (\varphi^{\ast}, w^{\ast}) \cong \R^{\omega}$.
\end{example}

\begin{example}\label{example_quotient}
Let $E$ be the lcs $l_{\infty} = \{(x_{n})\in \R^{\N}: \sup_{n}|x_{n}| < \infty \}$ equipped with the topology of pointwise convergence.
If we consider the canonical linear mapping $\pi$ of restriction from $C_p(\beta{\N})$ into $\R^{\N}$,
then the image of $\pi$ is exactly $E$. Since the mapping $\pi$ is not quotient, by this way we cannot decide whether $E$ is linearly Eberlein-Grothendieck.  
However, it has been proved that the space $C_p(\beta{\N})$ admits (another) linear continuous and quotient mapping onto $E$ \cite[Theorem 1]{BKS}.
Hence $E$ is linearly Eberlein-Grothendieck by Theorem \ref{prop_operations_2}.

Alternatively, we can argue as follows. Look at the set $K= \{\frac{1}{n}\delta_{n}: n\in\N\} \cup \{0\} \subset E^{\ast}$.
Then $K$ is compact in  $(E^{\ast}, w^{\ast})$ and $\span(K) = E^{\ast}$. Therefore, $E$ is linearly Eberlein-Grothendieck by Theorem \ref{theor:dual_compgen}.
\end{example}

\begin{example} \label{D(Omega)}
If $\Omega\subset \R^{n}$ is an open set, then the space of test functions $\mathfrak{D}(\Omega)$ is a complete Montel $(LF)$-space.
 As usual, $\mathfrak{D}'(\Omega)$ denotes its strong dual, the space of distributions. 
Not $\mathfrak{D}(\Omega)$ nor $\mathfrak{D}'(\Omega)$ is metrizable, they are not even sequential (see \cite {Dudley}). 
Applying Proposition \ref{prop_Montel} we conclude that $(E, w)$ is not Eberlein-Grothendieck for $E$ both $\mathfrak{D}(\Omega)$ and $\mathfrak{D}'(\Omega)$.
Also, by Corollary \ref{cor_original}, both $\mathfrak{D}(\Omega)$ and $\mathfrak{D}'(\Omega)$ considered with the original topologies are not Eberlein-Grothendieck.
\end{example}

\end{document}